\documentclass{amsart}

\usepackage{amsmath,amsthm,amssymb,amsfonts}

\theoremstyle{plain}
\newtheorem{theo+}           {Theorem}      [section]
\newtheorem{prop+}  [theo+]  {Proposition}
\newtheorem{lemm+}  [theo+]  {Lemma}
\newtheorem{cor+}  [theo+]  {Corollary}
\newenvironment{theorem}{\begin{theo+}}{\end{theo+}}

\newenvironment{corollary}{\begin{cor+}}{\end{cor+}}

\newcommand{\la}{\lambda}

\newcommand{\Ga}{\Gamma}
\newcommand{\De}{\Delta}
\begin{document}
\baselineskip 18pt
\larger[2]
\title[Karlsson--Minton type hypergeometric functions]
{Karlsson--Minton type hypergeometric\\ functions on the root system $C_n$}
\author{Hjalmar Rosengren}
\address{Department of Mathematics\\ Chalmers University of Technology and 
G\"oteborg University\\SE-412~96 G\"oteborg, Sweden}
\email{hjalmar@math.chalmers.se}

\begin{abstract}
We prove a reduction formula for Karlsson--Minton type hypergeometric
series on the root system $C_n$ and derive some consequences of this
identity. In particular, when combined with a  similar reduction
formula for $A_n$, it implies a $C_n$ Watson transformation due to
Milne and Lilly. 
\end{abstract}
\maketitle   

\section{Introduction}
The Karlsson--Minton summation formula \cite{mi,ka} 
is the hypergeometric identity
\begin{equation}\label{km}
{}_{p+2}F_{p+1}\left(\begin{matrix}a,b,c_1+m_1,\dots,c_p+m_p\\
b+1,c_1,\dots,c_p\end{matrix}\,;1\right)=\frac{\Ga(b+1)\Ga(1-a)}{\Ga(1+b-a)}
\prod_{i=1}^p\frac{(c_i-b)_{m_i}}{(c_i)_{m_i}},\end{equation}
which holds for $m_i$ non-negative integers and $\operatorname{Re}\,(a+|m|)<1$.
Accordingly, hyper\-geometric series with integral parameter differences have 
been called Karlsson--Minton type hypergeometric series; other results
for such series may be found in \cite{c,gk,g,s,se}.

In recent work \cite{r}, we have derived a very general reduction formula
for series of Karlsson--Minton type. We recall it in its most general form
as  \eqref{rg} below; here we only state a very
degenerate case, namely,
\begin{equation}\label{r}\begin{split}&
{}_{p+2}F_{p+1}\left(\begin{matrix}a,b,c_1+m_1,\dots,c_p+m_p\\
d,c_1,\dots,c_p\end{matrix}\,;1\right)
=\frac{\Ga(d)\Ga(d-a-b)}{\Ga(d-a)\Ga(d-b)}
\prod_{i=1}^p\frac{(c_i+1-d)_{m_i}}{(c_i)_{m_i}}\\
&\qquad
\times\frac{(a)_{|m|}}{(1+a+b-d)_{|m|}}
\sum_{x_1,\dots,x_p=0}^{m_1,\dots,m_p}\Bigg(\prod_{1\leq i<j\leq p}
\frac{c_i+x_i-c_j-x_j}{c_i-c_j}\frac{(b+1-d)_{|x|}}
{(1-|m|-a)_{|x|}}\\
&\qquad\times\prod_{i=1}^p\frac{(c_i-a)_{x_i}}{(1+c_i-d)_{x_i}}
\prod_{i,k=1}^p\frac{(c_i-c_k-m_k)_{x_i}}{(1+c_i-c_k)_{x_i}}\Bigg),
\end{split}\end{equation}
which holds for $m_i$ non-negative integers and  
$\operatorname{Re}\,(a+|m|+b-d)<0$. Note that the case $d=b+1$ is
\eqref{km}; similarly, the more general identity \eqref{rg}
implies a large number of results for Karlsson--Minton type hypergeometric
series from the papers mentioned above. 

The right-hand side
of \eqref{r} is a multivariable hypergeometric sum on the root
system $A_n$, a type  of series that were first introduced by
 Biedenharn, Holman and 
Louck \cite{hbl}, motivated by $6j$-symbols of the group $\mathrm{SU}(n)$. 
During  the last 25 years, hypergeometric series on
root systems has been a very
 active area of research 
with many applications. In the more general identity \eqref{rg}, both
the left- and the right-hand side are hypergeometric
series on $A_n$ (with different dimension $n$).
The connection between Karlsson--Minton type
series and hypergeometric series on $A_n$ encoded in \eqref{rg} 
allows one to recover many identities also for the latter type of series;
for instance, if we let $b=0$ in \eqref{r}, so that the left-hand side 
is $1$, we obtain the case $q=1$ of one of Milne's
 multivariable $q$-Saalsch\"utz summations \cite[Theorem 4.1]{m}.
More generally, \eqref{rg} implies 
 multivariable ${}_{10}W_9$  transformations due to Milne and Newcomb
\cite{mn} and  Kajihara \cite{k}.

The purpose of the present paper is to find results analogous to
those of \cite{r} for the root system $C_n$. In \cite{r}, our starting
point was an $A_n$ ${}_6\psi_6$ summation of Gustafson \cite{gu}; here we use
instead Gustafson's $C_n$  ${}_6\psi_6$ sum from \cite{gu2}.
Our main result is Theorem \ref{t}, which reduces a very general
multilateral Karlsson--Minton type series on $C_n$ to a finite sum.
In Section 4 we state some corollaries of Theorem \ref{t}. These include
a transformation and a summation formula for Karlsson--Minton
type series on $C_n$, Corollaries \ref{kmt} and \ref{kms}, respectively.
One special case of Theorem \ref{t} is a  rather curious
 transformation formula between finite sums,  Corollary \ref{ak}. 
Another interesting case is when the Karlsson--Minton type series
is one-variable, that is, connected to the root system $C_1$. In agreement
with the coincidence of root systems $C_1=A_1$, the same series may arise 
as a left-hand side of \eqref{rg}. This leads to a transformation formula
relating finite hypergeometric sums on the root systems $A_n$ and $C_n$,
which turns out to be a multivariable Watson transformation due to
Milne and Lilly \cite{ml}, given here as Corollary \ref{cml}.

\section{Preliminaries}
We will work with $q$-series rather than classical hypergeometric series,
with  $q$ a fixed complex number such that 
$0<|q|<1$. We will use the standard notation of \cite{gr}, but since
$q$ is fixed we suppress it from the notation. Thus we write
(this must not be confused with the standard notation for classical 
hypergeometric series used in the introduction)
$$(a)_k=\prod_{j=0}^\infty\frac{1-aq^j}{1-aq^{j+k}}
=\begin{cases}(1-a)(1-aq)\dotsm(1-aq^{k-1}),& k\geq 0,\\
\displaystyle\frac{ 1}
{(1-aq^{-1})(1-aq^{-2})\dotsm(1-aq^k)},& k<0,\end{cases}$$
$$(a_1,\dots,a_m)_k=(a_1)_k\dotsm(a_m)_k,$$
and analogously for infinite products $(a)_\infty=\prod_{j=0}^\infty(1-aq^j)$.

For $z=(z_1,\dots,z_n)\in\mathbb C^n$  we write $|z|=z_1+\dots+z_n$ and use the
corresponding capital letter to denote the product of the coordinates:
$Z=z_1\dotsm z_n$. 

We let
$$W(z)=\prod_{1\leq j\leq k\leq n}(1-z_jz_k)\prod_{1\leq j<k\leq n}(z_j-z_k),$$
which may be viewed as the Weyl denominator for the root system $C_n$
or the Lie group $\mathrm{Sp}(n)$. Basic hypergeometric series on $C_n$
are characterized by the factor
$$\frac{W(zq^y)}{W(z)}=\prod_{1\leq j\leq k\leq n}
\frac{1-z_jz_kq^{y_j+y_k}}{1-z_jz_k}
\prod_{1\leq j<k\leq n}\frac{q^{y_j}z_j-q^{y_k}z_k}{z_j-z_k}, $$
where the $z_j$ are fixed parameters and the $y_j$  summation indices. 
It will be useful to note that
\begin{equation}\label{ui}
\frac{W(z_1,\dots,z_{n-1},\la z_n)}{W(z_1,\dots,z_n)}
=\frac{1-\la^2z_n^2}{1-z_n^2}\prod_{k=1}^{n-1}\frac{(1-\la z_nz_k)(1-\la z_n/z_k)}
{(1-z_nz_k)(1-z_n/z_k)}.\end{equation}
We will also write
$$\frac{\De(zq^y)}{\De(z)}
=\prod_{1\leq j<k\leq n}\frac{q^{y_j}z_j-q^{y_k}z_k}{z_j-z_k}; $$
this factor characterizes hypergeometric functions on $A_n$.

For comparison we  recall the main result of \cite{r}, namely, the 
identity 
\begin{multline}\label{rg}
\begin{split}&
\sum_{\substack{y_1,\dots,y_n=-\infty\\y_1+\dots+y_n=0}}^\infty
\frac{\De(zq^y)}{\De(z)}
\prod_{\substack{1\leq k\leq n\\1\leq i\leq p}}\frac{(c_iz_kq^{m_i})_{y_k}}
{(c_iz_k)_{y_k}}\prod_{i,k=1}^n\frac{(a_iz_k)_{y_k}}
{(b_iz_k)_{y_k}}\\
&\quad=
\frac{(q^{1-|m|}/AZ,q^{1-n}BZ)_\infty}{(q,q^{1-|m|-n}B/A)_\infty}
\prod_{i,k=1}^n\frac{(b_i/a_k,qz_k/z_i)_\infty}{(q/a_kz_i,b_iz_k)_\infty}
\prod_{\substack{1\leq k\leq n\\1\leq i\leq p}}\frac{(q^{-m_i}b_k/c_i)_{m_i}}
{(q^{1-m_i}/c_iz_k)_{m_i}}\\
&\quad\quad\times
\sum_{x_1,\dots,x_p=0}^{m_1,\dots,m_p}\frac{\De(cq^x)}{\De(c)}\,q^{|x|}
\frac{(q^{n}/BZ)_{|x|}}{(q^{1-|m|}/AZ)_{|x|}}
\prod_{\substack{1\leq k\leq n\\1\leq i\leq p}}\frac{(c_i/a_k)_{x_i}}
{(qc_i/b_k)_{x_i}}\prod_{i,k=1}^p\frac{(q^{-m_k}c_i/c_k)_{x_i}}
{(qc_i/c_k)_{x_i}}.
\end{split}\end{multline}
Here, $m_i$ are non-negative integers, $|q^{1-|m|-n}B/A|<1$, and
it is assumed that no  denominators vanish.
The case $m_i\equiv 0$ of \eqref{rg} is Gustafson's
$A_n$ analogue of Bailey's ${}_6\psi_6$ summation \cite{gu}.
The proof of \eqref{rg} is based on induction on the $m_i$, with Gustafson's
identity as the starting point.

In the present paper we will imitate the analysis of \cite{r}, starting from
Gustafson's $C_n$ Bailey summation \cite{gu2},
which we write as
\begin{equation}\label{gi}\begin{split}&
\sum_{y_1,\dots,y_n=-\infty}^\infty\frac{W(zq^y)}{W(z)}
\prod_{\substack {1\leq k\leq n\\1\leq j\leq 2n+2}}
\frac{(a_jz_k)_{y_k}}{(qz_k/a_j)_{y_k}}\left(\frac qA\right)^{|y|}\\
&\qquad=\frac{\prod_{1\leq j\leq k\leq n}(qz_jz_k,q/z_jz_k)_\infty
\prod_{j,k=1}^n(qz_k/z_j)_\infty}
{\prod_{ {1\leq k\leq n,\,1\leq j\leq 2n+2}}(qz_k/a_j,q/a_jz_k)_\infty
}\frac{\prod_{1\leq j<k\leq 2n+2}(q/a_ja_k)_\infty}{(q/A)_\infty}.
\end{split}\end{equation}
This holds for $|q/A|<1$, as long as no denominators
 vanish. The case  $n=1$ of \eqref{gi}
 is Bailey's ${}_6\psi_6$ summation formula \cite[Equation (II.33)]{gr}.

\section{A reduction formula}

Our main result is the following identity, which reduces a very
general multilateral Karlsson--Minton type hypergeometric series
on the root system $C_n$ to a finite sum.

\begin{theorem}\label{t}
Let  $m_j$ be non-negative integers and
$a_j$, $c_j$,  $z_j$ parameters 
 such that $|q^{1-|m|}/A|<1$ and
none of the denominators in \eqref{te} vanishes.
Then the following identity holds:
\begin{multline}\label{te}
\sum_{y_1,\dots,y_n=-\infty}^\infty\frac{W(zq^y)}{W(z)}
\prod_{\substack{1\leq k\leq n\\1\leq j\leq p}}
\frac{(q^{m_j}c_jz_k,qz_k/c_j)_{y_k}}{(c_jz_k,q^{1-m_j}z_k/c_j)_{y_k}}
\prod_{\substack {1\leq k\leq n\\1\leq j\leq 2n+2}}
\frac{(a_jz_k)_{y_k}}{(qz_k/a_j)_{y_k}}\left(\frac{q^{1-|m|}}A\right)^{|y|}\\
\begin{split}&=\frac{\prod_{1\leq j\leq k\leq n}(qz_jz_k,q/z_jz_k)_\infty
\prod_{j,k=1}^n(qz_k/z_j)_\infty}
{\prod_{ {1\leq k\leq n,\,1\leq j\leq 2n+2}}(qz_k/a_j,q/a_jz_k)_\infty
}\frac{\prod_{1\leq j<k\leq 2n+2}(q/a_ja_k)_\infty}{(q/A)_\infty}\\
&\quad\times\frac{\prod_{{1\leq k\leq 2n+2,\,1\leq j\leq p}}
(c_ja_k)_{m_j}}
{\prod_{{1\leq k\leq n,\,1\leq j\leq p}}
(c_jz_k,c_j/z_k)_{m_j}}
\frac{\prod_{1\leq j<k\leq p}(c_jc_k)_{m_j+m_k}}
{\prod_{j,k=1}^p(c_jc_k)_{m_j}}\frac 1{(A)_{|m|}}\\
&\quad\times\sum_{x_1,\dots,x_p=0}^{m_1,\dots,m_p}
\frac{W(q^{-\frac 12}cq^x)}{W(q^{-\frac 12}c)}
\prod_{\substack {1\leq k\leq 2n+2\\1\leq j\leq p}}\frac{(c_j/a_k)_{x_j}}
{(c_ja_k)_{x_j}}\prod_{j,k=1}^p\frac{(q^{-1}c_jc_k,q^{-m_k}c_j/c_k)_{x_j}}
{(qc_j/c_k,q^{m_k}c_jc_k)_{x_j}}\left(Aq^{|m|}\right)^{|x|}.
\end{split}\end{multline}
\end{theorem}

The condition  $|q^{1-|m|}/A|<1$ ensures that that the series on the
left-hand side is absolutely convergent, so that the series manipulations
occurring in the proof are justified.

\begin{proof}
To organize the computations, it will be convenient to write, 
for $a\in \mathbb C^{2n+2}$, $z\in\mathbb C^n$, $y\in\mathbb Z^n$,
$c\in \mathbb C^p$ and
$m,\,x\in\mathbb N^p$,
$$P_y(a,z,c,m)=\frac{W(zq^y)}{W(z)}
\prod_{\substack{1\leq k\leq n\\1\leq j\leq p}}
\frac{(q^{m_j}c_jz_k,qz_k/c_j)_{y_k}}{(c_jz_k,q^{1-m_j}z_k/c_j)_{y_k}}
\prod_{\substack {1\leq k\leq n\\1\leq j\leq 2n+2}}
\frac{(a_jz_k)_{y_k}}{(qz_k/a_j)_{y_k}}\left(\frac{q^{1-|m|}}A\right)^{|y|},$$
$$U(a,z)=\frac{\prod_{1\leq j\leq k\leq n}(qz_jz_k,q/z_jz_k)_\infty
\prod_{j,k=1}^n(qz_k/z_j)_\infty}
{\prod_{ {1\leq k\leq n,\,1\leq j\leq 2n+2}}(qz_k/a_j,q/a_jz_k)_\infty
}\frac{\prod_{1\leq j<k\leq 2n+2}(q/a_ja_k)_\infty}{(q/A)_\infty}, $$
$$V(a,z,c,m)=\frac{\prod_{{1\leq k\leq 2n+2,\,1\leq j\leq p}}
(c_ja_k)_{m_j}}
{\prod_{{1\leq k\leq n,\,1\leq j\leq p}}
(c_jz_k,c_j/z_k)_{m_j}}
\frac{\prod_{1\leq j<k\leq p}(c_jc_k)_{m_j+m_k}}
{\prod_{j,k=1}^p(c_jc_k)_{m_j}}\frac 1{(A)_{|m|}}, $$
$$Q_x(a,c,m)=\frac{W(q^{-\frac 12}cq^x)}{W(q^{-\frac 12}c)}
\prod_{\substack {1\leq k\leq 2n+2\\1\leq j\leq p}}\frac{(c_j/a_k)_{x_j}}
{(c_ja_k)_{x_j}}\prod_{j,k=1}^p\frac{(q^{-1}c_jc_k,q^{-m_k}c_j/c_k)_{x_j}}
{(qc_j/c_k,q^{m_k}c_jc_k)_{x_j}}\left(Aq^{|m|}\right)^{|x|}, $$
so that \eqref{te} may be written as
\begin{equation}\label{tes}
\sum_{y_1,\dots,y_n=-\infty}^\infty P_y(a,z,c,m)=U(a,z)\,V(a,z,c,m)
\sum_{x_1,\dots,x_p=0}^{m_1,\dots,m_p}Q_x(a,c,m).\end{equation}

We will prove the theorem by induction on $|m|$, the case $m_j\equiv 0$
being Gustafson's identity \eqref{gi}. Thus we assume that \eqref{tes}
holds for fixed $p$ and $|m|$ and the other para\-meters free. Since both
sides are invariant under simultaneous permutations of the $m_j$ and
$c_j$, it is enough to prove that \eqref{tes} also holds when $m$ is
replaced by 
$$m^+=(m_1,\dots,m_{p-1},m_p+1).$$

It will be convenient to replace $n$ by $n+1$ in \eqref{tes} and write
\begin{gather*}a=(a_1,\dots,a_{2n+2}),\qquad 
z=(z_1,\dots,z_{n}),\qquad y=(y_1,\dots,y_{n}),\\
a^+=(a_1,\dots,a_{2n+4}),\qquad z^+=(z_1,\dots,z_{n+1}),\qquad
y^+=(y_1,\dots,y_{n+1}).\end{gather*}
We also specialize to the case $a_{2n+4}=z_{n+1}$.
Then the factor $1/(qz_{n+1}/a_{2n+4})_{y_{n+1}}$ 
on the left-hand side vanishes
unless $y_{n+1}\geq 0$, so that the series is supported on a half-space.
Next we let $a_{2n+3}= q/z_{n+1}$, which will cause most factors
involving $z_{n+1}$ to cancel. (In general it is not allowed
to put $a_{2n+3}= q/z_{n+1}$ in \eqref{te}, but the choice
of $a_{2n+4}$ removes this singularity.)  
We denote the corresponding left-hand side of \eqref{tes}
by $S$ and decompose it as
$$S=\sum_{y_1,\dots,y_{n+1}\in\mathbb Z,\,y_{n+1}\geq 0}
P_{y^+}(a^+,z^+,c,m)=\sum_{y_{n+1}=0}+\sum_{y_{n+1}\geq 1}=S_1+S_2.$$

Considering first $S_1$, we have
\begin{equation*}\begin{split}
\frac{W(z^+q^{(y,0)})}{W(z^+)}&=\frac{W(zq^y)}{W(z)}\prod_{k=1}^n
\frac{(1-q^{y_k}z_kz_{n+1})(1-q^{y_k}z_k/z_{n+1})}
{(1-z_kz_{n+1})(1-z_k/z_{n+1})}\\
&=\frac{W(zq^y)}{W(z)}\prod_{k=1}^n
\frac{(qz_kz_{n+1},qz_k/z_{n+1})_{y_k}}
{(z_kz_{n+1},z_k/z_{n+1})_{y_k}}.
\end{split}\end{equation*}
In particular, if we choose 
\begin{equation}\label{cz}z_{n+1}=q^{-m_p}/c_p,\end{equation}
the Weyl denominators combine with the factors involving $m_p$ as
$$\frac{W(z^+q^{(y,0)})}{W(z^+)}\prod_{k=1}^n
\frac{(q^{m_p}c_pz_k)_{y_k}}{(q^{1-m_p}z_k/c_p)_{y_k}}
=\frac{W(zq^y)}{W(z)}\prod_{k=1}^n
\frac{(q^{m_p+1}c_pz_k)_{y_k}}{(q^{-m_p}z_k/c_p)_{y_k}},$$
which gives
$$P_{(y,0)}(a^+,z^+,c,m)=P_y(a,z,c,m^+).$$
Thus, $S_1$ is a sum as in the theorem, but with $m$ replaced  by $m^+$.
To complete the induction, we must prove that
$$S_1=U(a,z)\,V(a,z,c,m^+)\sum_{x_1,\dots,x_p=0}^{m_1,\dots,m_{p-1},m_p+1}
Q_x(a,c,m^+).$$
This will be achieved by verifying the two identities
\begin{equation}\label{se}
S=U(a,z)\,V(a,z,c,m^+)\sum_{x_1,\dots,x_p=0}^{m_1,\dots,m_p}
Q_x(a,c,m^+),\end{equation}
\begin{equation}\label{se2}
S_2=-U(a,z)\,V(a,z,c,m^+)\sum_{x_1,\dots,x_{p-1}=0}^{m_1,\dots,m_{p-1}}
Q_{(x,m_p+1)}(a,c,m^+).\end{equation}

Starting with \eqref{se}, we know already that
$$S=U(a^+,z^+)\,V(a^+,z^+,c,m)\sum_{x_1,\dots,x_p=0}^{m_1,\dots,m_p}
Q_x(a^+,c,m). $$
It is not hard to check that (as explained above, when substituting
$a^+$ it is necessary to first let $a_{2n+4}=z_{n+1}$
and afterwards $a_{2n+3}=q/z_{n+1}$)
\begin{equation}\label{u}\begin{split}\frac{U(a^+,z^+)}{ U(a,z)}&
=\frac{1}{1-1/A}
\frac{\prod_{k=1}^{2n+2}(1-z_{n+1}/a_k)}
{\prod_{k=1}^{n+1}(1-z_{n+1}z_k)\prod_{k=1}^n(1-z_{n+1}/z_k)}\\
&=\frac{1}{1-A}
\frac{\prod_{k=1}^{2n+2}(1-a_k/z_{n+1})}
{\prod_{k=1}^{n+1}(1-1/z_{n+1}z_k)\prod_{k=1}^n(1-z_k/z_{n+1})},
\end{split}\end{equation}
\begin{equation}\label{v1}
\frac{V(a^+,z^+,c,m)}{ V(a,z,c,m)}=\frac{1-A}{1-Aq^{|m|}}
\prod_{j=1}^p\frac{1-q^{m_j}c_j/z_{n+1}}{1-c_j/z_{n+1}},\end{equation}
\begin{equation}\label{qa}\frac{Q_x(a^+,c,m)}{Q_x(a,c,m)}=
q^{|x|}\prod_{j=1}^p\frac{(1-c_j/z_{n+1})(1-q^{-1}c_jz_{n+1})}
{(1-q^{x_j}c_j/z_{n+1})(1-q^{x_j-1}c_jz_{n+1})}.\end{equation}
Combining this with the similarly derived identities
\begin{equation}\label{v2}\begin{split}\frac{V(a,z,c,m^+)}{V(a,z,c,m)}&
=\frac{1}{1-Aq^{|m|}}
\frac{\prod_{k=1}^{2n+2}(1-q^{m_p}c_pa_k)}
{\prod_{k=1}^n(1-q^{m_p}c_pz_k)(1-q^{m_p}c_p/z_k)}\\
&\quad\times\frac{\prod_{j=1}^{p-1}(1-q^{m_j+m_p}c_jc_p)}
{\prod_{j=1}^p(1-q^{m_p}c_jc_p)},\end{split}\end{equation}
$$\frac{Q_x(a,c,m^+)}{Q_x(a,c,m)}=
q^{|x|}\prod_{j=1}^p\frac{(1-q^{-m_p-1}c_j/c_p)(1-q^{m_p}c_jc_p)}
{(1-q^{x_j-m_p-1}c_j/c_p)(1-q^{x_j+m_p}c_jc_p)}$$
and using $z_{n+1}=q^{-m_p}/c_p$, we find that
$$U(a^+,z^+)\,V(a^+,z^+,c,m)\,Q_x(a^+,c,m)= U(a,z)\,V(a,z,c,m^+)\,
Q_x(a,c,m^+),$$
which proves \eqref{se}.

Next we show that $S_2$ is  a sum of the same type of $S$. 
The choice \eqref{cz} of $z_{n+1}$ corresponds to a removable singularity
of $S_2$. Namely, we must write
\begin{equation}\label{rs}
\frac{(q^{m_p}c_pz_{n+1})_{y_{n+1}}}{(c_pz_{n+1})_{y_{n+1}}}
=\frac{(q^{y_{n+1}}c_pz_{n+1})_{m_p}}{(c_pz_{n+1})_{m_p}}
=\frac{(q^{y_{n+1}-m_p})_{m_p}}{(q^{-m_p})_{m_p}},\end{equation}
which vanishes for $1\leq y_{n+1}\leq m_p$. To obtain a sum with 
$y_{n+1}\geq 0$ we  therefore replace $y_{n+1}$ with $y_{n+1}+m_p+1$
in the summation. This gives rise to factors of the form
$$(\lambda z_{n+1})_{y_k+m_p+1}\prod_{k=1}^{n}(\lambda z_k)_{y_k}=
(\lambda z_{n+1})_{m_p+1}\prod_{k=1}^{n+1}(\lambda w_k)_{y_k}, $$
where
$$w^+=(w_1,\dots,w_{n+1})=(z_1,\dots,z_n,q^{m_p+1}z_{n+1})
=(z_1,\dots,z_n,q/c_p).$$
Thus, the change of summation variables gives
$$S_2=\sum_{y_1,\dots,y_{n+1}\in\mathbb Z,\,y_{n+1}\geq 1}
P_{y^+}(a^+,z^+,c,m)
=M\sum_{y_1,\dots,y_n\in\mathbb Z,\,y_{n+1}\geq 0}P_{y^+}(a^+,w^+,c,m),$$
where, using \eqref{rs} with $y_{n+1}$ replaced by $1+m_p$,
\begin{equation}\label{m}\begin{split}M&=\frac{W(w^+)}{W(z^+)}
\frac{(q)_{m_p}}{(q^{-m_p})_{m_p}}\prod_{j=1}^{p-1}
\frac{(q^{m_j}c_jz_{n+1})_{m_p+1}}
{(c_jz_{n+1})_{m_p+1}}\prod_{j=1}^p
\frac{(qz_{n+1}/c_j)_{m_p+1}}
{(q^{1-m_j}z_{n+1}/c_j)_{m_p+1}}\\
&\quad\times
\prod_{k=1}^{2n+2}\frac{(a_kz_{n+1})_{m_p+1}}{(qz_{n+1}/a_k)_{m_p+1}}
\left(\frac{q^{-|m|}}{A}\right)^{m_p+1}.\end{split}\end{equation}

We first  rewrite the multiplier $M$. By \eqref{ui},
$$\frac{W(w^+)}{W(z^+)}=\frac{1-q^2/c_p^2}{1-q^{-2m_p}/c_p^2}
\prod_{k=1}^n\frac{(1-qz_k/c_p)(1-q/z_kc_p)}{(1-q^{-m_p}z_k/c_p)
(1-q^{-m_p}/z_kc_p)}.$$
Plugging this into \eqref{m}, and using  \eqref{cz} and
the standard identities 
$$\frac{(a)_n}{(b)_n}=\left(\frac ab\right)^n
\frac{(q^{1-n}/a)_n}{(q^{1-n}/b)_n},
\qquad \frac{(q)_n}{(q^{-n})_n}=(-1)^nq^{n(n+1)/2}, $$
we obtain
\begin{equation}\label{c1}\begin{split}
M&=(-1)^{m_p}\left(Aq^{|m|-\frac 12 m_p}\right)^{m_p+1}
\frac{1-q^{-2}c_p^2}{1-q^{2m_p}c_p^2}
\prod_{k=1}^n \frac{(1-q^{-1}c_pz_k)(1-q^{-1}c_p/z_k)}{(1-q^{m_p}c_pz_k)
(1-q^{m_p}c_p/z_k)}\\
&\quad\times
\prod_{k=1}^{2n+2}\frac{(c_p/a_k)_{m_p+1}}{(q^{-1}c_pa_k)_{m_p+1}}
\prod_{j=1}^{p-1}\frac{(q^{-m_j}c_p/c_j)_{m_p+1}}{(c_p/c_j)_{m_{p}+1}}
\prod_{j=1}^p\frac{(q^{-1}c_jc_p)_{m_p+1}}{(q^{m_j-1}c_jc_p)_{m_p+1}}.
\end{split}\end{equation}

By  our induction hypothesis, we have
\begin{equation}\label{c2}
S_2=M\, U(a^+,z^+)\,V(a^+,z^+,c,m)\sum_{x_1,\dots,x_p=0}^{m_1,\dots,m_p}
Q_x(a^+,c,m),\end{equation}
where $a_{2n+3}$ and $a_{2n+4}$ are related to $z_{n+1}$ as above, but where
instead of \eqref{cz} we have that
$z_{n+1}=q/c_p.$
Using \eqref{u}, \eqref{v1} and \eqref{v2} with $z_{n+1}=q/c_p$ gives
\begin{equation}\label{c3}\begin{split}&\frac{U(a^+,z^+)\,V(a^+,z^+,c,m)}
{U(a,z)\,V(a,z,c,m^+)}
=\frac{1-q^{2m_p}c_p^2}{1-q^{-2}c_p^2}
\prod_{k=1}^{2n+2}\frac{1-q^{-1}a_kc_p}{1-q^{m_p}a_kc_p}\\
&\qquad\times\prod_{k=1}^n
\frac{(1-q^{m_p}c_pz_k)(1-q^{m_p}c_p/z_k)}{(1-q^{-1}c_pz_k)(1-q^{-1}c_p/z_k)}
\prod_{j=1}^p\frac{(1-q^{m_j-1}c_jc_p)(1-q^{m_p}c_jc_p)}
{(1-q^{m_j+m_p}c_jc_p)(1-q^{-1}c_jc_p)}.\end{split}\end{equation}
When $z_{n+1}=q/c_p$, \eqref{qa} vanishes unless $x_p=0$, in which case
\begin{equation}\label{c4}\frac{Q_{(x,0)}(a^+,c,m)}{Q_{(x,0)}(a,c,m)}=q^{|x|}
\prod_{j=1}^{p-1}\frac{(1-c_j/c_p)(1-q^{-1}c_jc_p)}
{(1-q^{x_j}c_j/c_p)(1-q^{x_j-1}c_jc_p)}.\end{equation}
Finally we want to compare $Q_{(x,0)}(a,c,m)$ and
$Q_{(x,m_p+1)}(a,c,m^+)$. Again using \eqref{ui}, we have
$$\frac{W(q^{-\frac 12}cq^{(x,m_{p}+1)})}{W(q^{-\frac 12}cq^{(x,0)})}
=\frac{1-q^{2m_p+1}c_p^2}{1-q^{-1}c_p^2}\prod_{j=1}^{p-1}
\frac{(1-q^{x_j+m_p}c_jc_p)(1-q^{1-x_j+m_p}c_p/c_j)}
{(1-q^{x_j-1}c_jc_p)(1-q^{-x_j}c_p/c_j)}, $$
which gives, after simplifications,
\begin{multline}\label{c5}\begin{split}&
\frac{Q_{(x,m_p+1)}(a,c,m^+)}{Q_{(x,0)}(a,c,m)}
=(-1)^{m_p+1}q^{|x|}\left(Aq^{|m|-\frac 12 m_p}\right)^{m_p+1}
\frac{(c_p^2)_{m_p}}{(q^{m_p+1}c_p^2)_{m_p}}\\
&\times
\prod_{k=1}^{2n+2}\frac{(c_p/a_k)_{m_p+1}}{(c_pa_k)_{m_p+1}}
\prod_{j=1}^{p-1}\left(\frac{(1-q^{-1}c_jc_p)(1-c_j/c_p)}
{(1-q^{x_j-1}c_jc_p)(1-q^{x_j}c_j/c_p)}
\frac{(c_jc_p,q^{-m_j}c_p/c_j)_{m_p+1}}{(q^{m_j}c_jc_p,c_p/c_j)_{m_p+1}}
\right).\end{split}\end{multline}
Combining the equations \eqref{c1}, \eqref{c2}, \eqref{c3}, \eqref{c4}
and \eqref{c5}, we eventually obtain \eqref{se2}. This completes the proof.
\end{proof}

\section{Corollaries}

In this section we point out some interesting consequences and special cases
of Theorem \ref{t}. Throughout, it is assumed that the $m_j$ are non-negative
integers and that no denominators  vanish.

One of the most conspicuous features of 
\eqref{te} is that the sum on the right is independent of the parameters
$z_j$. This  implies a transformation formula
for the series on the left. The case $n=1$ of the resulting identity is
due to Schlosser \cite[Corollary 8.6]{se}. Schlosser also gave a 
generalization to the root system $A_n$  \cite[Theorem 4.2]{s}; 
cf.~also \cite[Corollary 4.2]{r}.

\begin{corollary}\label{kmt}
For $|q^{1-|m|}/A|<1$, the following identity holds:
\begin{multline*}\begin{split}&
\sum_{y_1,\dots,y_n=-\infty}^\infty\frac{W(zq^y)}{W(z)}
\prod_{\substack{1\leq k\leq n\\1\leq j\leq p}}
\frac{(q^{m_j}c_jz_k,qz_k/c_j)_{y_k}}{(c_jz_k,q^{1-m_j}z_k/c_j)_{y_k}}
\prod_{\substack {1\leq k\leq n\\1\leq j\leq 2n+2}}
\frac{(a_jz_k)_{y_k}}{(qz_k/a_j)_{y_k}}\left(\frac{q^{1-|m|}}A\right)^{|y|}\\
&\qquad=\prod_{1\leq j\leq k\leq n}\frac{(qz_jz_k,q/z_jz_k)_\infty}
{(qw_jw_k,q/w_jw_k)_\infty}\prod_{j,k=1}^n\frac{(qz_k/z_j)_\infty}
{(qw_k/w_j)_\infty}\\
&\qquad\quad\times\prod_{\substack{1\leq k\leq n\\1\leq j\leq 2n+2}}
\frac{(qw_k/a_j,q/a_jw_k)_\infty}{(qz_k/a_j,q/a_jz_k)_\infty}
\prod_{\substack{1\leq k\leq n\\1\leq j\leq p}}
\frac{(c_jw_k,c_j/w_k)_{m_j}}{(c_jz_k,c_j/z_k)_{m_j}}\end{split}\\
\times\sum_{y_1,\dots,y_n=-\infty}^\infty\frac{W(wq^y)}{W(w)}
\prod_{\substack{1\leq k\leq n\\1\leq j\leq p}}
\frac{(q^{m_j}c_jw_k,qw_k/c_j)_{y_k}}{(c_jw_k,q^{1-m_j}w_k/c_j)_{y_k}}
\prod_{\substack {1\leq k\leq n\\1\leq j\leq 2n+2}}
\frac{(a_jw_k)_{y_k}}{(qw_k/a_j)_{y_k}}\left(\frac{q^{1-|m|}}A\right)^{|y|}.
\end{multline*}
\end{corollary}

If we assume that  $a_{n+j}=a_j^{-1}$ for $1\leq j\leq n$
and choose $w_j=a_j$ in Corollary \ref{kmt}, 
the factor $(w_{n+j}a_j)_{y_j}/(qw_j/a_j)_{y_j}
=(1)_{y_j}/(q)_{y_j}$ on the right vanishes for $y_j\neq 0$, so that
the sum reduces to $1$.  Alternatively, we  may in this situation use
  Corollary \ref{n0} below to compute the
right-hand side of \eqref{te}. Writing
 $a_{2n+1}=b$, $a_{2n+2}=d$, either of these methods
 gives the following identity.
When $n=1$, it reduces to an identity of  Chu \cite{c}.
For $A_n$ analogues of Chu's identity, cf.~\cite[Corollary 4.3]{s},
\cite[Corollaries 4.3 and 4.4]{r}.

\begin{corollary}\label{kms}
For $|q^{1-|m|}/bd|<1$, the following identity holds:
\begin{multline*}
\sum_{y_1,\dots,y_n=-\infty}^\infty\Bigg(\frac{W(zq^y)}{W(z)}
\prod_{\substack{1\leq k\leq n\\1\leq j\leq p}}
\frac{(q^{m_j}c_jz_k,qz_k/c_j)_{y_k}}{(c_jz_k,q^{1-m_j}z_k/c_j)_{y_k}}
\prod_{j,k=1}^n
\frac{(a_jz_k,z_k/a_j)_{y_k}}{(qa_jz_k,qz_k/a_j)_{y_k}}\\
\begin{split}&\quad\times\prod_{k=1}^n
\frac{(bz_k,dz_k)_{y_k}}{(qz_k/b,qz_k/d)_{y_k}}
\left(\frac{q^{1-|m|}}{bd}\right)^{|y|}\Bigg)\\
&=
\prod_{j,k=1}^n\frac{(qz_k/z_j,qa_k/a_j)_\infty}{(qz_k/a_j,qz_ka_j,
q/a_jz_k,qa_j/z_k)_\infty}
\prod_{k=1}^n\frac{(qa_k/b,q/a_kb,qa_k/d,q/a_kd)_\infty}
{(qz_k/b,q/z_kb,qz_k/d,q/z_kd)_\infty}\\
&\quad\times\prod_{1\leq j\leq k\leq n}
(qz_jz_k,q/z_jz_k)_\infty\prod_{1\leq j< k\leq n}
(qa_ja_k,q/a_ja_k)_\infty
\prod_{\substack{1\leq k\leq n\\1\leq j\leq p}}\frac{(c_ja_k,c_j/a_k)_{m_j}}
{(c_jz_k,c_j/z_k)_{m_j}}.
\end{split}\end{multline*}
\end{corollary}

To obtain an identity closer to the original Karlsson--Minton
summation formula \eqref{km} one should specialize the parameters in Corollary
\ref{kms} so that the summation indices  are bounded from
below. Essentially, this forces $n=2$, when we may choose 
$b=z_1$, $d=z_2$. The resulting identity seems interesting
enough to write out explicitly; it is a $C_2$ version of 
Gasper's well-poised Karlsson--Minton type summation from \cite{g}.

\begin{corollary}
For $|q^{1-|m|}/z_1z_2|<1$, the following identity holds:
\begin{multline*}
\sum_{y_1,y_2=0}^\infty\frac{W(zq^y)}{W(z)}
\prod_{\substack{1\leq k\leq 2\\1\leq j\leq p}}
\frac{(q^{m_j}c_jz_k,qz_k/c_j)_{y_k}}{(c_jz_k,q^{1-m_j}z_k/c_j)_{y_k}}
\prod_{j,k=1}^2
\frac{(a_jz_k,z_k/a_j,z_jz_k)_{y_k}}{(qa_jz_k,qz_k/a_j,qz_k/z_j)_{y_k}}
\left(\frac{q^{1-|m|}}{z_1z_2}\right)^{|y|}\\
=\frac{(qz_1^2,qz_1z_2,qz_2^2,qa_1a_2,q/a_1a_2)_\infty}{(q/z_1z_2)_\infty}
\prod_{j,k=1}^2\frac{(qa_k/a_j)_\infty}{(qz_k/a_j,qz_ka_j)_\infty}
\prod_{\substack{1\leq k\leq 2\\1\leq j\leq p}}\frac{(c_ja_k,c_j/a_k)_{m_j}}
{(c_jz_k,c_j/z_k)_{m_j}}.
\end{multline*}
\end{corollary}

More generally, one may
 choose the parameters so that the summation indices
 on the  left-hand side of \eqref{te} are bounded from
below or above. A particularly symmetric case arises when
both these conditions hold, so that we have a  finite sum.
To this end we choose $a_j=z_j$, $a_{n+j}=q^{-l_j}/z_j$, $1\leq j\leq n$
 in Theorem \ref{t}
and write  $a_{2n+1}=b$, $a_{2n+2}=d$. 
Since we have a rational identity in $b$,
$d$ the condition $|q^{1-|m|}/A|=|q^{1-|m|+|l|}/bd|<1$ is then
superfluous. The resulting identity is reminiscent of
 transformation formulas for $A_n$ hypergeometric series recently obtained
by Kajihara \cite{k}.

 \begin{corollary}\label{ak}
The following identity holds:
\begin{multline*}
\sum_{y_1,\dots,y_n=0}^{l_1,\dots,l_n}\Bigg(\frac{W(zq^y)}{W(z)}
\prod_{\substack{1\leq k\leq n\\1\leq j\leq p}}
\frac{(q^{m_j}c_jz_k,qz_k/c_j)_{y_k}}{(c_jz_k,q^{1-m_j}z_k/c_j)_{y_k}}
\prod_{j,k=1}^n
\frac{(z_jz_k,q^{-l_j}z_k/z_j)_{y_k}}{(qz_k/z_j,q^{1+l_j}z_jz_k)_{y_k}}\\
\begin{split}&\quad\times
\prod_{k=1}^n\frac{(bz_k,dz_k)_{y_k}}{(qz_k/b,qz_k/d)_{y_k}}
\left(\frac{q^{1-|m|+|l|}}{bd}\right)^{|y|}\Bigg)\\
&=\frac{(q/bd)_{|l|}}{(q^{-|l|}bd)_{|m|}}
\frac{\prod_{j=1}^p(c_jb,c_jd)_{m_j}}{\prod_{k=1}^n(qz_k/b,qz_k/d)_{l_k}}
\frac{\prod_{j,k=1}^n(qz_jz_k)_{l_k}}
{\prod_{1\leq j<k\leq n}(qz_jz_k)_{l_j+l_k}}
\frac{\prod_{1\leq j<k\leq p}(c_jc_k)_{m_j+m_k}}
{\prod_{j,k=1}^p(c_jc_k)_{m_j}}\\
&\quad\times
\prod_{\substack{1\leq k\leq n\\1\leq j\leq p}}\frac{(q^{-l_k}c_j/z_k)_{m_j}}
{(c_j/z_k)_{m_j}}
\sum_{x_1,\dots,x_p=0}^{m_1,\dots,m_p}\Bigg(
\frac{W(q^{-\frac 12}cq^x)}{W(q^{-\frac 12}c)}
\prod_{\substack {1\leq k\leq n\\1\leq j\leq p}}
\frac{(q^{l_k}c_jz_k,c_j/z_k)_{x_j}}
{(c_jz_k,q^{-l_k}c_j/z_k)_{x_j}}\\
&\quad\times\prod_{j,k=1}^p\frac{(q^{-1}c_jc_k,q^{-m_k}c_j/c_k)_{x_j}}
{(qc_j/c_k,q^{m_k}c_jc_k)_{x_j}}\prod_{j=1}^p
\frac{(c_j/b,c_j/d)_{x_j}}{(c_jb,c_jd)_{x_j}}
\left(q^{|m|-|l|}bd\right)^{|x|}\Bigg).
\end{split}
\end{multline*}
\end{corollary}

Theorem \ref{t} has some interesting consequences for low values of
$n$. An inspection of the proof shows that it
 holds for $n=0$, if
the left-hand side of \eqref{te} is interpreted as $1$.
This leads to a $C_n$ analogue of the terminating ${}_6W_5$ summation
formula, which is in fact  the special
case $a_j=z_j$, $a_{n+j}=q^{-m_j}/z_j$, $1\leq j\leq n$ of Gustafson's 
identity \eqref{gi} (cf.~also \cite{lm}), or, equivalently, the case
$m_j\equiv 0$ of Corollary \ref{ak}. We include it here as 
 a first illustration of how  Theorem~\ref{t} is related to known
results for ``classical'' (i.e.~not of Karlsson--Minton type)
 $C_n$ hypergeometric series. Compared to Theorem \ref{t}
we have replaced $p$ with $n$ and $c_j$, $a_1$, $a_2$ with $q^{1/2}z_j$, 
$q^{1/2}a$, $q^{1/2}b$, respectively.

\begin{corollary}[Gustafson]\label{n0}
The following identity holds:
\begin{equation*}\begin{split}&
\sum_{x_1,\dots,x_n=0}^{m_1,\dots,m_n}
\frac{W(zq^x)}{W(z)}
\prod_{j=1}^n\frac{(z_j/a,z_j/b)_{x_j}}
{(qz_ja,qz_jb)_{x_j}}\prod_{j,k=1}^n\frac{(z_jz_k,q^{-m_k}z_j/z_k)_{x_j}}
{(qz_j/z_k,q^{m_k+1}z_jz_k)_{x_j}}\left(abq^{|m|+1}\right)^{|x|}\\
&\qquad=\frac{(qab)_{|m|}}{\prod_{j=1}^n(qaz_j,qbz_j)_{m_j}}
\frac{\prod_{j,k=1}^n(qz_jz_k)_{m_j}}
{\prod_{1\leq j<k\leq n}(qz_jz_k)_{m_j+m_k}}.\end{split}\end{equation*}
\end{corollary}

Next we consider the case $n=1$ of Theorem \ref{t}, 
when the left-hand side of \eqref{te}
is a one-variable well-poised ${}_{2p+6}\psi_{2p+6}$ series. 
After letting $z_1=a^{1/2}$,
$(a_1,a_2,a_3,a_4)=a^{-1/2}(b,c,d,e)$, $c_j=qa^{1/2}/f_j$, 
it takes the following form.

\begin{corollary}\label{fl}
For $|a^2q^{1-|m|}/bcde|<1$, the following identity holds:
\begin{multline}\label{ci}\begin{split}
&\sum_{y=-\infty}^\infty
\frac{1-aq^{2y}}{1-a}\frac{(b,c,d,e)_y}{(aq/b,aq/c,aq/d,aq/e)_y}
\prod_{j=1}^p\frac{(f_j,aq^{1+m_j}/f_j)_{y}}{(q^{-m_j}f_j,aq/f_j)_{y}}
\left(\frac{a^{2}q^{1-|m|}}{bcde}\right)^{y}\\
&\quad=\frac{(q,aq,q/a,aq/bc,aq/bd,aq/be,aq/cd,aq/ce,aq/de)_\infty}
{(q/b,q/c,q/d,q/e,aq/b,aq/c,aq/d,aq/e,a^{2}q/bcde)_\infty}
\frac{1}{(bcde/a^2)_{|m|}}\\
&\quad\quad\times\prod_{j=1}^p\frac{(qb/f_j,qc/f_j,qd/f_j,qe/f_j)_{m_j}}
{(aq/f_j,q/f_j)_{m_j}}\frac
{\prod_{1\leq j<k\leq p}(aq^2/f_jf_k)_{m_j+m_k}}
{\prod_{j,k=1}^p(aq^2/f_jf_k)_{m_j}}\\
&\quad\quad\times\sum_{x_1,\dots,x_p=0}^{m_1,\dots,m_p}
\Bigg(\frac{W(\sqrt{aq}\,q^x/f)}{W(\sqrt{aq}/f)}
\prod_{j=1}^p\frac{(aq/bf_j,aq/cf_j,aq/df_j,aq/ef_j)_{x_j}}
{(bq/f_j,cq/f_j,dq/f_j,eq/f_j)_{x_j}}\\
&\quad\quad\times\prod_{j,k=1}^p
\frac{(aq/f_jf_k,q^{-m_k}f_k/f_j)_{x_j}}{(qf_k/f_j,aq^{2+m_k}/f_jf_k)_{x_j}}
\left(\frac{bcdeq^{|m|}}{a^2}\right)^{|x|}\Bigg).
\end{split}\end{multline}
\end{corollary}
We remark that the factor
$$\frac{W(\sqrt{aq}\,q^x/f)}{W(\sqrt{aq}/f)}=\prod_{1\leq j\leq k\leq p}
\frac{1-aq^{x_j+x_k+1}/f_jf_k}{1-aq/f_jf_k}\prod_{1\leq j< k\leq p}
\frac{q^{x_j}/f_j-q^{x_k}/f_k}{1/f_j-1/f_k}$$
does not depend on the choice of square root.

Corollary \ref{fl} may be compared with Corollary 4.11 of \cite{r}, 
which is just the case $n=2$ of \eqref{rg}. 
It says that the left-hand side of \eqref{ci} equals
\begin{multline*}\begin{split}
&\frac{(q,aq,q/a,aq/bc,aq/bd,aq/be,aq/cd,aq/ce,aq/de)_\infty}
{(q/b,q/c,q/d,q/e,aq/b,aq/c,aq/d,aq/e,a^{2}q/bcde)_\infty}\\
&\times\prod_{j=1}^p
\frac{(bq/f_j,cq/f_j)_{m_j}}
{(aq/f_j,q/f_j)_{m_j}}\frac{(de/a)_{|m|}}{(bcde/a^2)_{|m|}}\\
&\times\sum_{x_1,\dots,x_p=0}^{m_1,\dots,m_p}
\frac{\De(q^x/f)}{\De(1/f)}
\frac{(bc/a)_{|x|}}{(aq^{1-|m|}/de)_{|x|}}
\prod_{j=1}^p\frac{(aq/df_j,aq/ef_j)_{x_j}}
{(qb/f_j,qc/f_j)_{x_j}}\prod_{j,k=1}^p\frac{(q^{-m_k}f_k/f_j)_{x_j}}
{(qf_k/f_j)_{x_j}}\,q^{|x|}.\end{split}\end{multline*}
That this quantity equals the
right-hand side of \eqref{ci} is equivalent to
a multivariable Watson transformation due to
Milne and Lilly \cite[Theorem 6.6]{ml}. 
After replacing $p$ with $n$, $f_j$ with
$\sqrt{aq}/z_j$ and $(b,c,d,e)$ with $\sqrt{aq}\,
(b^{-1},c^{-1},d^{-1}, e^{-1})$, it takes the following form.

\begin{corollary}[Milne and Lilly]\label{cml}
One has the identity
\begin{multline*}\begin{split}&\sum_{x_1,\dots,x_n=0}^{m_1,\dots,m_n}
\Bigg(\frac{W(zq^x)}{W(z)}
\prod_{j=1}^n\frac{(bz_j,cz_j,dz_j,ez_j)_{x_j}}
{(qz_j/b,qz_j/c,qz_j/d,qz_j/e)_{x_j}}\\
&\quad\quad\times\prod_{j,k=1}^n
\frac{(z_jz_k,q^{-m_k}z_j/z_k)_{x_j}}{(qz_j/z_k,q^{1+m_k}z_jz_k)_{x_j}}
\left(\frac{q^{|m|+2}}{bcde}\right)^{|x|}\Bigg)\\
&\quad=\frac{(q/de)_{|m|}\prod_{j,k=1}^n(qz_jz_k)_{m_j}}
{\prod_{j=1}^n(qz_j/d,qz_j/e)_{m_j}\prod_{1\leq j<k\leq n}(qz_jz_k)_{m_j}}
\\
&\quad\quad\times\sum_{x_1,\dots,x_n=0}^{m_1,\dots,m_n}
\frac{\De(zq^x)}{\De(z)}
\frac{(q/bc)_{|x|}}{(q^{-|m|}de)_{|x|}}\prod_{j=1}^n\frac{(dz_j,ez_j)_{x_j}}
{(qz_j/b,qz_j/c)_{x_j}}\prod_{j,k=1}^n\frac{(q^{-m_k}z_j/z_k)_{x_j}}
{(qz_j/z_k)_{x_j}}\,q^{|x|}.
\end{split}\end{multline*}
\end{corollary}

The  proof of Corollary \ref{cml}
 obtained here  gives a nice explanation of 
 why  such a transformation formula,  relating a $C_n$ ${}_8W_7$ series
and an $A_n$
${}_4\phi_3$, exists: on the level of 
Karlsson--Minton type hypergeometric series it  reflects the coincidence of 
root systems $A_1=C_1$. It is
appropriate to remark here that our proof of 
Theorem \ref{t} depended on guessing the explicit expression for the
right-hand side of \eqref{te}. 
This task was much simplified by having access to
the Milne--Lilly transformation, and thus 
(given also the results of \cite{r}) knowing  the identity
 in advance for $n=1$.

A generalization of Corollary \ref{cml} to the level of multivariable 
balanced ${}_{10}W_9$ series has been obtained by Bhatnagar and Schlosser 
\cite[Theorem 2.1]{bs}. We have not been able to obtain this identity in 
our approach.
 Note that, when $n=2$, the right-hand side of \eqref{te}
is a $p$-variable ${}_{10}W_9$, but to make it balanced one must let
$A=q^{1-|m|}$, which corresponds to a pole of the left-hand side.


\begin{thebibliography}{99} 
\bibitem[BS]{bs} G.\ Bhatnagar and M.\ Schlosser,
{\it $C_n$ and $D_n$ very-well-poised ${}_{10}\phi_9$ transformations},
Constr.\ Approx.\ 14 (1998),  531--567.

\bibitem[C]{c} W.\ Chu,  
{\it Partial-fraction expansions and well-poised bilateral series},
Acta Sci.\ Math.\ (Szeged) 64 (1998),  495--513. 



\bibitem[G1]{gk} G.\ Gasper,  
{\it Summation formulas for basic hypergeometric series}, 
SIAM J.\ Math.\ Anal.\ 12 (1981),  196--200. 

\bibitem[G2]{g}
G.\ Gasper,  
{\it Elementary derivations of summation and transformation formulas for 
$q$-series}, 
Special Functions, $q$-Series and Related Topics, 55--70, 
Fields Inst.\ Commun.\ 14,  Providence, RI, 1997. 

\bibitem[GR]{gr} G.\ Gasper and M.\ Rahman, Basic Hypergeometric Series,
Cambridge University Press, Cambridge, 1990. 

\bibitem[Gu1]{gu} R.\ A.\ Gustafson, {\it Multilateral summation theorems 
for ordinary and basic hypergeometric series in ${\rm U}(n)$}, 
SIAM J.\ Math.\ Anal.\ 18 (1987),  1576--1596. 


\bibitem[Gu2]{gu2} R.\ A.\ Gustafson,
{\it The Macdonald identities for affine root systems of classical type and 
hypergeometric series very-well-poised on
semisimple Lie algebras}, 
Ramanujan International Symposium on Analysis (Pune, 1987), 185--224, 
Macmillan of India, New Delhi, 1989. 

\bibitem[HBL]{hbl} W.\ J.\ Holman, L.\ C.\ Biedenharn and J.\ D.\ Louck,  
{\it On hypergeometric series well-poised in ${\rm SU}(n)$}, 
SIAM J.~Math.\ Anal.\ 7 (1976),  529--541. 


\bibitem[K]{k} Y.\ Kajihara, {\it Euler transformation formulas for multiple
basic hypergeometric series of type $A$ and some applications},
Adv.~Math., to appear.



\bibitem[Ka]{ka} P.\ W.\ Karlsson, {\it Hypergeometric functions with 
integral parameter differences},
J.\ Math.\ Phys.\ 12 (1971), 270--271. 

\bibitem[LM]{lm} G.\ M. Lilly and S.\ C.\ Milne, 
{\it The $C_l$ Bailey transform and Bailey lemma},
Constr.\ Approx.\ 9 (1993),  473--500. 

\bibitem[M]{m}S.\ C.\ Milne, {\it Balanced ${}_3\phi_2$ summation theorems 
for ${\rm U}(n)$ basic hypergeometric series},
Adv.\ Math.\ 131 (1997),  93--187. 

\bibitem[ML]{ml}  S.\ C.\ Milne and G.\ M. Lilly, 
{\it Consequences of the $A_l$ and $C_l$ Bailey transform and
Bailey lemma},  Discrete Math.\ 139 (1995), 319--346. 

\bibitem[MN]{mn} S.\ C.\ Milne and J.\ W.\ Newcomb, 
{\it  ${\rm U}(n)$ very-well-poised ${}_{10}\phi_9$ transformations},
J. Comput. Appl. Math. 68 (1996), 239--285.


\bibitem[Mi]{mi} B.\ Minton, {\it 
Generalized hypergeometric function of unit argument},
J.\ Math.\ Phys.\ 11 (1970), 1375--1376. 

\bibitem[R]{r} H. Rosengren, 
{\it Reduction formulae for Karlsson--Minton type hypergeometric functions}, 
\verb+math.CA/0202232+. 

\bibitem[S1]{s}M.\ Schlosser, {\it Multilateral transformations of $q$-series 
with quotients of parameters that are nonnegative integral powers of $q$}, 
 Contemp.\ Math.\ 291 (2001), 203--227. 

\bibitem[S2]{se} M.\ Schlosser, {\it Elementary derivations of identities 
for bilateral basic hypergeometric series}, Selecta Math., to appear.


\end{thebibliography}
\end{document}